\documentclass[leqno,12pt]{article}

\usepackage{inputenc}

\usepackage[english]{babel}

\usepackage{amsmath}

\usepackage{amssymb}

\usepackage{amsthm}

\usepackage{textcomp}

\usepackage[all]{xy}

\date{}

\newtheorem{lem}{Lemma}
\newtheorem{cor}[lem]{Corollary}
\newtheorem{prop}[lem]{Proposition}
\newtheorem{thm}[lem]{Theorem}

\theoremstyle{definition}

\newtheorem{rem}[lem]{Remark}
\newtheorem{ex}[lem]{Example}

\numberwithin{lem}{section}

\usepackage{soul,color}

\title{On the Hilbert function of the tangent cone of a monomial curve}

\newcommand{\ap}{\mathrm {Ap}}

\newcommand{\ord}{\mathrm {ord}}

\newcommand{\gr}{\mathrm{gr}}

\newcommand{\m}{\mathfrak{m}}

\newcommand{\N}{\mathbb{N}}

\newcommand{\Z}{\mathbb{Z}}

\newcommand{\Apery}{Ap\'ery}

\author{M. D'Anna \footnotemark[1]\thanks{Universit\`a di Catania, Dipartimento di Matematica e
Informatica, Viale A. Doria, 6, 95125 Catania, Italy} \and M. Di
Marca \footnotemark[2]\thanks{Scuola Superiore di Catania, Via
Valdisavoia, 9, 95123 Catania, Italy} \and V. Micale
\footnotemark[1]}

\begin{document}

\maketitle

\begin{abstract}

\medskip

\noindent In this paper we study the Hilbert function of
$\gr_{\mathfrak{m}}(R)$, when $R$ is a numerical semigroup ring
or, equivalently, the coordinate ring of a monomial curve. In
particular, we prove a sufficient condition for a numerical
semigroup ring in order get a non-decreasing Hilbert function,
without making any assumption on its embedding dimension;
moreover, we show how this new condition allows to improve known
results about this problem. To this aim we use certain invariants
of the semigroup, with particular regard to its \Apery-set.


\end{abstract}

\section*{Introduction}
Let $(R,\mathfrak{m})$ be a Noetherian local ring with
$|R/\mathfrak{m}|=\infty$ and let
$\gr_{\mathfrak{m}}(R)=\oplus_{i\geq 0}\mathfrak{m}^i/
\mathfrak{m}^{i+1}$ be the associated graded ring of $R$ with
respect to $\mathfrak{m}$. The study of the properties of
$\gr_{\mathfrak{m}}(R)$ is a classical subject in local algebra,
not only in the general $d$-dimensional case, but also under
particular hypotheses (that allow to obtain more precise results).
A classical problem in this context is to study the Hilbert
function of $R$, i.e., by definition, the Hilbert function of
$\gr_\m(R)$.

\medskip

In this paper we are interested in the Hilbert function of $R$,
when $R$ is a numerical semigroup ring. The study of numerical
semigroup rings is motivated by their connection to singularities
of monomial curves and by the possibility of translating algebraic
properties into numerical properties (see e.g. \cite{BDF}).
However, even in this particular case, many pathologies occur,
hence these rings are also a great source of interesting examples.

From the geometrical point of view, given a numerical semigroup
$S$ generated by $n$ coprime integers $g_1,g_2,\ldots,g_n$, the
numerical semigroup ring, $R=k[[S]]$ is the completion of the
local ring at the origin of the monomial curve
$C=C(g_1,\ldots,g_n)$ parameterized by $x_1={t^{g_1}}, \dots
x_n=t^{g_n}$. Hence its associated graded ring is the coordinate
ring of the tangent cone of $C$ in the origin. Moreover,
$\gr_\m(R)$ is isomorphic to the ring $k[x_1,x_2, \ldots,
x_n]/I(C)_*$, where $I(C)$ is the defining ideal of $C$ and
$I(C)_*$ is the ideal generated by the homogeneous terms of least
degree of the polynomials in $I(C)$.

\medskip

One classical problem about the Hilbert function is to find
conditions on $R$ or $\gr_\m(R)$ to get a non decreasing function.
In the context of one-dimensional local rings, it is well known
that the Hilbert function is not decreasing if the embedding
dimension is at most $3$ (see \cite{E} and \cite{EB}) and
counterexamples (for reduced one-dimensional rings) were given for
embedding dimension bigger than or equal to $4$ (see \cite{O} and
\cite{GR}). However, for semigroup rings, there are no examples of
decreasing Hilbert function when the embedding dimension is
smaller than $10$; so the problem is still open for semigroup
rings $R$ with embedding dimension $4\leq e.d.(R) \leq 9$.

Another open problem in the one-dimensional case was posed by
Rossi in \cite{Ro}: if $R$ is a Gorenstein one-dimensional local
ring, is it true that the Hilbert function of $\gr_\m(R)$ is not
decreasing? In the context of semigroup rings, it is equivalent to
ask whether the Hilbert function of $k[[S]]$, with $S$ symmetric,
is non-decreasing.

The problem if the Hilbert function of a semigroup ring is
non-decreasing, has been extensively studied. If $\gr_\m(R)$ is
Cohen-Macaulay, then the pro\-blem becomes trivial thanks to a
result of A. Garcia in \cite{Ga}. For the general case, recent
results can be found, e.g., in \cite{A-M-S}, where many families
of non-decreasing Hilbert function of semigroup rings are obtained by using the
technique of gluing semigroups (see also \cite{JZ}), and in
\cite{P-T}, where the authors study particular $4$-generated
semigroups. Furthermore, in \cite{CJZ} new results on this problem
are obtained introducing the \Apery-table of a semigroup; in
particular, the authors proved that if $S$ is $4$-generated and if
the tangent cone is Buchasbaum, then $k[[S]]$ has non-decreasing
Hilbert function.

The main result of this paper is a sufficient condition for a
numerical semigroup ring in order to get a non-decreasing Hilbert
function, without making any assumption on its embedding dimension
(see Theorem \ref{main} and Corollary \ref{mainb}); to this aim we
use certain invariants of the semigroup, with particular regard to
its \Apery-set. Successively, a careful use of the proof of the
main result allows us to get a computationally more efficient,
necessary condition for the decreasing of the Hilbert function
(see Proposition \ref{prop} and the subsequent remark). Finally,
we show how these results can be applied to improve known results
about this problem (see Corollaries \ref{cor0}, \ref{cor2} and
\ref{ed45}).

\section{Preliminaries }

We start this section recalling some well known facts on numerical
semigroups and semigroup rings. For more details see, e.g.,
\cite{BDF}.

A \textit{numerical semigroup} $S$ is a subsemigroup of
$(\mathbb{N},+)$ that includes $0$. There is a natural partial
order on $S$ that is defined as follows: let $a,b \in S$, then
\[a \leq_S b \ \Longleftrightarrow \ \exists \ u \in S\ :\ a+u=b.\]
\noindent The set of minimal elements with respect to this order
is called \textit{minimal set of generators} of $S$. It is always
finite because, by definition, for any $s \in S, s \neq 0$, two
minimal generators have to be different modulo $s$. Once fixed the
minimal set of generators, each element of $S$ can be written as
finite sum of these elements. Hence $S$ is determined by its
minimal set of generators. We denote by $\langle g_1, g_2, \ldots,
g_n\rangle$ the numerical semigroup $S$ whose minimal set of
generators is $\{g_1, g_2, \ldots, g_n\}$, where $g_1 < g_2<
\ldots < g_n$. Since the semigroup $S=\langle g_1, g_2, \ldots,
g_n \rangle$ is isomorphic to $\langle dg_1, dg_2, \ldots, dg_n
\rangle$ for any $d \in \mathbb{N}\setminus 0$, we can assume that
$\gcd(g_1, g_2, \ldots, g_n)=1$. This is equivalent to say that
$|\N \setminus S|$ is finite, so it is well defined the maximum of
the numbers that does not belong to the semigroup, called
\textit{Frobenius number of $S$} and denoted by $f$.

From now on, we will call a numerical semigroup simply semigroup.

\medskip

A \textit{relative ideal} of a semigroup $S$ is a set $H \subset
\mathbb Z$, $H \neq \emptyset$, such that $H+S \subseteq H$ and
$H+s \subseteq S$, for some $s \in S$; if $H \subseteq S$, it is
called \textit{ideal}. If $H$ and $L$ are relative ideals, then
also $kH=\{h_1+h_2+\ldots+h_k : h_1,h_2,\ldots,h_k \in H\}$ (for
all $k \in \N$) and $H-_{\Z}L=\{z \in \Z : z+l \in H,\ \forall l
\in L\}$ are relative ideals. The ideal $M=S \setminus \{0\}$ is
called \textit{maximal ideal}.

\medskip Let $k$ be an infinite field
and let $S=\langle g_1, g_2, \ldots, g_n \rangle$; the ring
$R=k[[t^S]]=k[[t^{g_1}$,$t^{g_2},\ldots,t^{g_n}]]$ is called
\textit{semigroup ring associated to $S$}. The ring $R$ is a
one-dimensional local domain, with maximal ideal
$\m=(t^{g_1},t^{g_2} \ldots, t^{g_n})$ and quotient field
$k((t))$. Considering the $\m$-adic filtration, let $gr^h(R)$ be
the quotient $\m^h/\m^{h+1}$. From the direct sum of the $gr^h(R)$
we obtain the \textit{associated graded ring} $\gr_\m(R)$
explicitly defined as $\gr_\m(R)=\bigoplus_{h\geq
0}\m^h/\m^{h+1}$. Setting $k=R/\m$, the \textit{Hilbert function}
of $R$ is then given by
\[H_R(h)=\dim_k gr^h(R), \ \ \forall \ n \in \N.\]

There exists a strong connection between a semigroup and its
associated ring. In fact, through the \textit{natural valuation}
function $v:k((t))\rightarrow \Z \cup \infty$, that is
\[v\big(\sum_{n=i}^{\infty}r_nt^n\big)=i, \ \ i \in \Z, \ \ r_i\neq 0,\]
we get $v(R)=S$ and many other properties. For example, if $I$ and
$J$ are fractional ideal of $R$, then $v(I)$ and $v(J)$ are
relative ideal of $S$ and so are $v(I \cap J)$, $v(I:J)$ and
$v(I^n)$ for all $n \in \N$. Furthermore, if $I$ and $J$ are
monomial fractional ideals, the following relations hold: $v(I
\cap J)=v(I)\cap v(J)$, $v(I:J)=v(I)-v(J)$ and $v(I^n)=nv(I)$.
Moreover, if $J \subseteq I$ are fractional ideals of $R$, then:
\[\dim_k(I/J)=|v(I)\setminus v(J)|\ .\]These facts hold, in particular, for $I=\m^h$ and $J=\m^{h+1}$,
for all $h \geq 0$. Therefore, since $v(\m)=M$, the Hilbert
function of the semigroup ring $H_R$ is equivalent to the
\textit{Hilbert function of S} which is
\[H_S(h)=|hM \setminus (h+1)M|\ , \ \ \forall \ h\in \N\](when $h=0$ we set, as usual, $\m^0=R$ and $0M=S$).

\medskip

We denote by $\ap(S)$ the \Apery-set of $S$ with respect to the
smallest generator $g_1$, which is the set $\{\omega_0, \omega_1,
\ldots, \omega_{g_1-1}\}$, where $\omega_i=min\{s \in S : s\equiv
i \ (\text{mod} \ g_1)\}$.

The numerical semigroup $S'=\bigcup_i(iM-_{\Z}iM)$ is called
\textit{blow up of $S$} and it corresponds to the blow up of $R$.
By \cite[Lemma 1]{N}, $S$ is generated by $\{g_1, g_2-g_1, \ldots,
g_n-g_1\}$ (but this is not necessarily its minimal set of
generators). Moreover $S' \supseteq hM-hg_1=\{s-hg_1 \ : \ s \in
hM\}$ (for every $h\geq 1$) and the equality holds for every $h$
large enough. The \Apery-set of $S'$ with respect to $g_1$ is
denoted by $\ap(S')=\{\omega'_0, \omega'_1, \ldots,
\omega'_{g_1-1}\}$.

We recall two important sets of invariants of $S$, introduced by
Barucci and Fr\"oberg in \cite{BF}. For each $i=0,1, \ldots,
g_1-1$, let $a_i$ be the only integer such that
$\omega'_i+a_ig_1=\omega_i$ and let $b_i=\max\{l: \omega_i \in
lM\}$. Clearly $b_0=a_0=0$. Furthermore, Barucci and Fr\"oberg
proved that $1 \leq b_i \leq a_i$ for all $i$ (see \cite{BF},
Lemma 2.4) and that $\gr_\m(R)$ is Cohen-Macaulay if and only if
the equality $a_i=b_i$ holds for each $i$ (see \cite{BF}, Theorem
2.6).

We will need also to consider another set of invariants introduced
in \cite{DMM} and related to the previous ones: $c_i=\min\{h\ :\
\omega'_i \in hM-hg_1\}$. We have that $a_i \leq c_i$ (for every
$i$) and that $b_i<a_i$ if and only if $a_i<c_i$ (see \cite{DMM},
Proposition 3.5).

\section{Relation between \Apery-set and Hilbert \\ function}

In this section we relate the coefficient $a_i$ and $b_i$
introduced in the previous section with the value of the Hilbert
function.

From the characterization of $H_R$ given in the first section, it
is obvious that $H_R$ is non-decreasing if and only if
$$
|(h-1)M \setminus hM| \leq |hM \setminus (h+1)M|\ , \ \ \forall \
h \geq 1.$$

In order to study this inequality, it is natural to consider the
elements $s$ belonging to $(h-1)M \setminus hM$ and to add $g_1$:
now, if $s+g_1 \in hM \setminus (h+1)M$, for all $s \in (h-1)M
\setminus hM$, we get an injective function between the two sets.
If this is the case for every $h \geq 1$, then $H_R$ is
non-decreasing. We will see that this situation corresponds to the
case $\gr_\m(R)$ Cohen-Macaulay. On the contrary, it can happen
that the set

\[D_h:= \{s \in (h-1)M \setminus hM : s+g_1 \in (h+1)M\}\]

\noindent is non empty for some $h\ge 2$.

Let $r$ be the reduction number of $\m$, that is the minimal
natural number such that $\mathfrak m^{r+1}=x\mathfrak m^{r}$, with
$x$ a superficial element of $R$ (recall that such number $r$
exists by \cite[Theorem 1, Section 2]{NR}). We notice that, in the
semigroup ring case, the valuation of $x$ is necessarily
$v(x)=g_1$; hence the multiplicity of $R$, i.e.
$\dim_k(\m^h/\m^{h+1})=|hM \setminus (h+1)M|$ (for any $h\geq r$),
coincides with $g_1$. \\

Let $s \in S$; the maximal index $h$ such that $s \in hM$ is
called the \textit{order} of $s$ and it is denoted by $\ord(s)$.
We note that $\ord(\omega_i)=b_i$. Since $\ord(s)=h$ if and only
if s $\in hM \setminus (h+1)M$, we often say in this case that $s$
is \textit{on the h-th level}. We also say that an element $s$
\textit{skips the level when adding} $g_1$ if $\ord(s)=h$ and
$\ord(s+g_1)>h+1$. With this terminology, we can say that $D_h$ is
given by the elements on the $(h-1)$-th level that skip the level
when adding $g_1$.

Let $D=\cup_{i\ge 2} D_i$. Notice that the condition $D=\emptyset$
is equivalent to say that the image of $t^{g_1}$ in $\gr_\m(R)$ is
not a zero-divisor, i.e., by \cite{Ga}, $\gr_\m(R)$ is
Cohen-Macaulay.

The following result shows that the condition $a_i>b_i$ for
$\omega_i$ in $\ap(S)$ is related with the elements in $D$.

\begin{prop}\label{car}
Let $S$ be a semigroup. For each index $i$ there exists an element
$s \in D$, $s \equiv i\ (mod\ g_1)$ if and only if $a_i>b_i$. In
particular, if $a_i=b_i$ for every $i$, then $D=\emptyset$.
\end{prop}
\proof

$(\Rightarrow$) Let $s\in (h-1)M\setminus hM$ and let $s \equiv i \ (\text{mod} \
g_1)$; hence $s=\omega_i+\lambda g_1$. We use the
induction on $\lambda$ to prove that $a_i=b_i$ implies $s \notin D_h.$\\
For $\lambda=0$, we have $\omega_i \in b_iM \setminus (b_i+1)M$;
if we suppose that $\omega_i+g_1 \in (b_i+2)M$ we get
$\omega_i+g_1=\omega'_i+a_ig_1+g_1=(\omega'_i-g_1)+(b_i+2)g_1 \in
(b_i+2)M$,
that implies $\omega'_i-g_1 \in S'$, against the minimality of $\omega'_i$ in $S'.$\\
We now suppose the thesis true for $\lambda-1 \geq 0$ and we fix
$s=\omega_i+\lambda g_1$. By the inductive hypothesis $s\in
(b_i+\lambda)M \setminus (b_i+\lambda+1)M$; again, if we assume
$s+g_1 \in (b_i+\lambda+2)M$ then we could write
$s+g_1=(\omega'_i-g_1)+(b_i+\lambda+2)g_1 \in (b_i+\lambda+2)M$,
in contrast with the definition of $\omega'_i.$\\
$(\Leftarrow$) We know that $a_i>b_i$ implies $c_i>a_i.$ Let us
consider the element
$s'=\omega'_i+c_ig_1=\omega'_i+a_ig_1+(c_i-a_i)g_1$. By definition
of $c_i$, we have $s' \in c_iM$. Moreover,
$\omega_i=\omega'_i+a_ig_1 \in b_iM \setminus (b_i+1)M$ and
$s'=\omega_i+(c_i-a_i)g_1$. Hence the inequality $b_i+c_i-a_i<c_i$
implies that there must be an element $s=\omega_i+\lambda g_1$
(for some $\lambda$, $0 \leq \lambda \leq c_i-a_i-1$) that skips
the level when adding $g_1$, i.e. $s \in D$.
\endproof

In \cite[Theorem 2.6]{BF}, the authors prove (in the more general
context of one-dimensional analytically irreducible rings) that
$\gr_\m(R)$ is Cohen-Macaulay if and only if \ $a_i=b_i$, for
every $i=0,\dots,g_1-1$. As we noticed above, the
cohen-macaulaynness of $\gr_\m(R)$ is equivalent to $D=\emptyset$;
hence the previous proposition give a different proof of the
result of Barucci and Froberg in the semigroup ring case.
Moreover, as a corollary, we get the well known result that, if
$\gr_\m(R)$ is Cohen-Macaulay, then $H_R$ is non-decreasing.

\medskip
In order to compare $|(h-1)M \setminus hM|$ and $|hM \setminus
(h+1)M|$, we want to take into account the number of the skipping
elements in each level. Hence, for each $h\ge 1$, we set

\[C_h=\{s \in hM \setminus (h+1)M:\ s-g_1 \notin (h-1)M\setminus hM,\};\]
in other words, $C_h$ is the set of elements on the $h$-th level
which don't come from any element on the previous level by adding
$g_1$. This means that if $s\in C_h$, then either $s-g_1$ is a
skipping element coming from a level lower than the $(h-1)$-th
level, or $s$ is an element of $Ap(S)$ of order $h$.

We notice that the sets $C_h$ and $D_h$ arise naturally in this
context and were already defined in \cite{CJZ}; hence we conformed
our terminology to the names appearing in that paper.\\

With this notation it is straightforward that $|(h-1)M \setminus
hM| \leq |hM \setminus (h+1)M|$ if the number of elements on the
$(h-1)$-th level that skip level when adding $g_1$ is smaller than
or equal to the number of elements in $C_h$, i.e. $|D_h| \leq
|C_h|$. Hence $$H_R \ \ \text{is non-decreasing} \ \
\Longleftrightarrow \ \ |D_h|\leq|C_h|\ , \ \ \forall \ h \in
\{2,\dots, r\}
$$
(where $r$ is the reduction number of $\m$; notice also that for
$h=1$ it is always true that $|S\setminus M|=1 \leq |M\setminus 2M|$
or, equivalently, $D_1=\emptyset$).\\

Our next goal is to find conditions for determining an injective
function from $D_h$ to $C_h$. We recall that each element $s$ in
the semigroup can be written as linear combination of the
generators with coefficient in $\mathbb{N}$. If we have

\[s=\lambda_1g_1+\lambda_2g_2+\cdots+\lambda_ng_n, \quad \lambda_i\in \mathbb{N},\]

\noindent we say that this is a \textit{maximal representation} of
$s$ if $\sum_{i=1}^{n} \lambda_i=\ord(s)$. The maximal
representation of an element is not unique in general. If we have
two maximal representations of the same element $s$, we will write

\[\lambda_1g_1+\lambda_2g_2+\cdots+\lambda_ng_n \prec \lambda'_1g_1+\lambda'_2g_2+\cdots+\lambda'_ng_n\]

\noindent if $(\lambda_1,\lambda_2,\dots,\lambda_n)$ is smaller
than $(\lambda'_1,\lambda'_2,\dots,\lambda'_n)$ in the usual
lexicographic order in $\mathbb{N}^n$.

The following lemma is similar to \cite[Lemma 4.2(2)]{CJZ}. We
prove it in this form, since we will need it for the proof of the
main theorem.

\begin{lem}\label{lem}

For every index $h\ge 2$ there exists a function $\psi:D_h
\rightarrow C_h$.

\end{lem}

\proof

Let $s \in D_h$. Then $\ord(s)=h-1$ and $s+g_1 \in (h+1)M$; hence
there exists $k \geq 2$ such that $s+g_1 \in (h+k-1)M \setminus
(h+k)M$. Let

\[s+g_1=g_{l_1}+g_{l_2}+\cdots+g_{l_h}+g_{l_{h+1}}+\cdots+g_{l_{h+k-1}},\]

\noindent (with $g_{l_1}\leq g_{l_2}\leq \cdots \leq
g_{l_{h+k-1}}$) be the greatest among all the maximal
representations of $s+g_1$ with respect to the Lex order.

We define $\psi(s):=g_{l_1}+g_{l_2}+\cdots+g_{l_h}$; hence we have
$\psi(s) \in hM \setminus (h+1)M$, because it is part of a maximal
representation. Let $s':=\psi(s)-g_1$ and let us assume, by
contradiction, that $s' \in (h-1)M$. We get
$s+g_1=s'+g_1+g_{l_{h+1}}+\cdots+g_{l_{h+k-1}}$ that implies $s
=s'+g_{l_{h+1}}+\cdots+g_{l_{h+k-1}}\in hM$, against the fact that
$\ord(s)=h-1$. Hence $\psi(s)-g_1 \notin (h-1)M$ and $\psi(s) \in
C_h$.
\endproof

Now we give a sufficient condition for $S$ in order to have
$|(h-1)M \setminus hM| \leq |hM \setminus (h+1)M|$. The Example \ref{es} below will
illustrate the procedure of the proof of the next theorem.

\begin{thm}\label{main}

If $|D_h|\leq h+1$, then there exists an injective function
$\widetilde{\psi}:D_h \rightarrow C_h.$
\end{thm}

\proof

Let $s$ be an element in $D_h$, and let

\[s+g_1=g_{l_1}+g_{l_2}+\cdots+g_{l_h}+g_{l_{h+1}}+\cdots+g_{l_{h+k-1}},\]

\noindent (with $g_{l_1}\leq g_{l_2}\leq \cdots \leq
g_{l_{h+k-1}}$) be the greatest among all the maximal
representations of $s+g_1$ with respect to the Lex order.

As in Lemma \ref{lem} we map $s$ to
$\psi(s)=g_{l_1}+g_{l_2}+\cdots+g_{l_h}$. Let $J=|D_h|$ and let

\[\psi(D_h)=\{\psi(s_1), \psi(s_2), \ldots, \psi(s_J)\}\]

\noindent without removing the possible repetitions. Hence
$|\psi(D_h)|=|D_h| \leq h+1$. Let us order $\psi(D_h)$ according
to the decreasing Lex order:

\[\psi(D_h)=\{\psi_1, \psi_2, \ldots, \psi_J\}, \quad \psi_1\succeq \psi_2\succeq \cdots \succeq \psi_J,\]

\noindent where $\psi_j=\psi(s_m)$ for some $1 \leq m \leq J$.

If all the images $\psi(s_m)$ are pairwise different, we get the
thesis with $\widetilde{\psi}=\psi$. Otherwise, let $a$ be the
minimum among the indexes such that $\psi_a=\psi_{a+1}$ and let
$s_u$ and $s_v$ be the pre-images of $\psi_a$ and $\psi_{a+1}$,
respectively.

Let $\psi_a=\psi_{a+1}=g_{l_1}+g_{l_2}+\cdots+g_{l_h}$. Since $s_u
\neq s_v$, there exists a generator $g_p$, with $g_p
> g_{l_h}$, that appears in the maximal representations of
$s_u+g_1$ (or $s_v+g_1$) and does not appear in
$\psi_a=\psi_{a+1}$ (otherwise, if both $s_u+g_1$ and $s_v+g_1$
have maximal representations involving only $g_{l_1}, \dots
g_{l_h}$, by $\psi_a=\psi_{a+1}$ we would get $s_u+g_1 \leq _S
s_v+g_1$, or viceversa, and consequently $s_u \leq _S s_v$, or
viceversa; contradiction against the fact that $s_u$ and $s_v$
have both order $h-1$). Without loss of generality we can assume
that $g_p$ appears in the representations of $s_v+g_1$.

Now we can define a new function $\psi'$ on $D_h$ so that
$\psi'(s)=\psi(s)$, for every $s\in D_h$, $s \neq s_v$, and
$\psi'(s_v)=\psi_{a+1}-g_{l_h}+g_p=g_{l_1}+g_{l_2}+\cdots+g_{l_{h-1}}+g_p$.
\\

We now have $\psi'(s_v)\prec\psi(s_u)=\psi_a$. The new set of
images is

\[\psi'(D_h)=[\psi(D_h)\cup \{\psi'(s_v)\}]\setminus \{\psi_{a+1}\}.\]

\noindent We reorder $\psi'(D_h)$ and we rename the elements
$\psi'_{j}$, for $j=a+1,\ldots,J$ according to the decreasing Lex
order. We have:

\[\psi'_1\succ\psi'_2\succ \cdots \succ\psi'_a\succeq \psi'_{a+1}\succeq \cdots \succeq \psi'_J. \]

\noindent Again, if all the elements in $\psi'(D_i)$ are pairwise
different, we get the thesis. Otherwise, we repeat the same
argument as above by taking the minimum of the indices for which
we have an equality in the chain and we redefine the correspondent
images. We observe that this index could be $a$ again. In this
case, we are sure that $\psi_{a+1} \neq \psi'_{a+1}$ and we can
compare the two pre-images of $\psi'_{a}$ and $\psi'_{a+1}$ as in
the previous step. By redefining one of them, we get a new set of
images $\psi''(D_h)$.

There is the possibility that by ordering $\psi''(D_h)$ we find
again an equality for the index $a$. We note that this event can
happen at most $J-a$ times and no conditions are required to
redefine the function at each step.

After a finite number of steps (say $w \geq 1$) we will have:
\[\psi^{(w)}_1\succ\psi^{(w)}_2\succ \cdots \succ\psi^{(w)}_a\succ \psi^{(w)}_{a+1}\succeq
 \psi^{(w)}_{a+2}\cdots \succeq \psi^{(w)}_J. \]

\noindent When this condition is satisfied we say that we have
performed the first block of steps.

Again, if the images of $\psi^{(w)}$ are all pairwise different,
we get the thesis with $\widetilde{\psi}=\psi^{(w)}$. Otherwise,
let $b$ be the minimum among the indexes such that
$\psi^{(w)}_b=\psi^{(w)}_{b+1}$ (with $b>a$). Let again $s_u \neq
s_v$ be the pre-images of $\psi^{(w)}_b$ and $\psi^{(w)}_{b+1}$,
respectively; moreover, we set
$\psi_b^{(w)}=\psi_{b+1}^{(w)}=g_{l_1}+g_{l_2}+\cdots+g_{l_{h-1}}+g_{l_q}$
(with $g_{l_1}\leq g_{l_2}\leq \cdots \leq g_{l_{i-1}}\leq
g_{l_q}$). This time, since $s_u \neq s_v$, there exists a
generator $g_p$\ , with $g_p
> g_{l_{h-1}}$, that appears in the maximal representations of
$s_u+g_1$ (or $s_v+g_1$) (otherwise, if both $s_u+g_1$ and
$s_v+g_1$ have maximal representations involving only $g_{l_1},
\dots, g_{l_{h-1}}$, since $g_{l_1}\leq g_{l_2}\leq \cdots \leq
g_{l_{h-1}}$ are the first $h-1$ elements in the maximal
representation of both $s_u+g_1$ and  $s_v+g_1$, we would get
$s_u+g_1 \leq _S s_v+g_1$, or viceversa, and consequently $s_u
\leq _S s_v$, or viceversa; contradiction against the fact that
$s_u$ and $s_v$ have both order $h-1$. Notice that in this case we
are not able any more to compare $g_p$ and $g_{l_q}$). Without
loss of generality we can assume that $g_p$ appears in the
representations of $s_v+g_1$.

We can define a new function $\psi^{(w+1)}$ on $D_h$ so that
$\psi^{(w+1)}(s)=\psi^{(w)}(s)$, for every $s\in D_h$, $s \neq
s_v$, and
$\psi^{(w+1)}(s_v)=\psi^{(w)}_{b+1}-g_{l_{h-1}}+g_p=g_{l_1}+g_{l_2}+\cdots+g_{l_{h-2}}+g_p+g_{l_q}$
(or $g_{l_1}+g_{l_2}+\cdots+g_{l_{h-2}}+g_{l_q}+g_p$, if $g_p >
g_{l_q}$). This means that at this and at all the subsequent steps
we will rearrange the summands in non-decreasing order. As in the
previous block of steps, we go on until we get a $\psi^{(z)}$ such
that
\[\psi^{(z)}_1\succ\psi^{(z)}_2 \succ \cdots \succ\psi^{(z)}_a\succ \cdots \succ \psi^{(z)}_{b}
\succ\psi^{(z)}_{b+1} \succeq \psi^{(z)}_{b+2}\cdots \succeq
\psi^{(z)}_J. \]

\noindent When this condition is satisfied we say that we have
performed the second block of steps.

We would like to proceed until we obtain a chain of proper
inequalities; i.e., deno\-ting the last defined function by
$\widetilde{\psi}$, until we get:
$$\widetilde{\psi_1}\succ\widetilde{\psi_2}\succ\cdots\succ \widetilde{\psi_J}.$$

To this aim, in the worst case, we need to perform $J-a$ blocks of
steps, where $a$ is the index of the first equality. Since at the
$j$-th block of steps we substitute $g_{l_{h-j+1}}$ with the new
generator $g_p$, we are sure that we can perform $h$ blocks of
steps. Hence, in order to get the desired injective function, it
is sufficient that $J-a \leq h$. Since $a \geq 1$ and $J=|D_h|$ we
get the thesis.

We finally observe that every $\widetilde{\psi_a}$ is still an
element of $C_h$, coming from some $s\in D_h$. Furthermore, if
$\widetilde{\psi_a}\succ\widetilde{\psi_b}$, then
$\widetilde{\psi_a}\ne\widetilde{\psi_b}$, since we are assuming
that the summands of $\widetilde{\psi_a}$ (for every index $a$)
are in non-decreasing order.
Hence
$\widetilde{\psi}$ is an injective function.
\endproof

\begin{cor}\label{mainb}
If $|D_h|\leq h+1$ for every $h\ge 2$, then $H_R$ is
non-decreasing.
\end{cor}

The next example is appropriate to illustrate the procedure of the
proof of the main theorem.

\begin{ex}\label{es}
Let $S=\langle 24,25,36,51,54\rangle$. Its Hilbert function is
non-decreasing; in fact, it assumes the following values
$$1,5,11,16,19,20,21,22,22,22,22,23,24, \rightarrow.$$ We have
$|D_2|=1$, $|D_3|=3$, $|D_4|=4$, $|D_5|=4$ and $|D_h| \leq 3$ for
every $h \geq 5$. Hence it is fulfilled the condition of the
previous theorem and corollary. Let us analyze $D_5=\{s \in 4M
\setminus 5M : s+24 \in 6M\}=\{126,137,155,166\}$. We have
$$\aligned
126+24&=6\cdot 25 \\
137+24&=5\cdot 25+36 \\
155+24&=5\cdot 25+54 \\
166+24&=4 \cdot 25+36+54
\endaligned
$$
The function $\psi$ defined in Lemma \ref{lem}, gives $$\psi(126)=
\psi(137)= \psi(155)=5\cdot 25 \succ \psi(166)=4\cdot 25+36.$$
Following the proof of Theorem \ref{main}, we have $a=1$ and we
define $\psi'(137)=4\cdot25+36$: now we have $\psi'(126)=
\psi'(155)\succ \psi'(137)= \psi'(166)$. So again we have an
equality for $a=1$; we are forced to define $\psi''(155)= 4 \cdot
25 +54$ and we get: $\psi''(126)\succ
\psi''(137)=4\cdot25+36=\psi''(166)\succ \psi''(155)$. Now we have
completed the first block of steps and we have an equality for
$b=2$.

In the maximal representation of $166+24$, it appears $54$; hence
we can define $\psi^{(3)}(166)=3\cdot25+54+36$. After reordering
its summands we obtain $\psi^{(3)}(126)=5\cdot
25\succ\psi^{(3)}(137)=4\cdot25+36 \succ\psi^{(3)}(155)=4 \cdot 25
+54\succ\psi^{(3)}(166)=3\cdot 25+36+54$. Now we have completed
the second block of steps and, since we have inequalities for
every index of the chain, we can set $\widetilde \psi=\psi^{(3)}$.

Notice that $C_4=\{125,136,154,165,191\}=\widetilde \psi(D_4) \cup
\{191\}$ and that $191-24 \notin S$.
\end{ex}

\begin{rem}
We note that the result of Theorem \ref{main} is the best possible
as the semigroup $S=\langle 13,19,24,44,49,54,55,59,60,66 \rangle$
has $|D_2|=4$ and $H_R$ decreasing (see \cite{MT}). More precisely
we have:
$$\aligned
44+13&=19+19+19 \\
49+13&=19+19+24 \\
54+13&=19+24+24 \\
59+13&=24+24+24
\endaligned
$$
The function $\psi$, defined in Lemma \ref{lem}, gives $\psi(44)=
\psi(49)=19+19 \succ \psi(54)=19+24 \succ \psi(59)=24+24$. If we
try to follow the procedure of Theorem \ref{main}, we have to
define $\psi'(49)=19+24$ (in this case the first block of steps
consists of one step); now we have the equality
$\psi'(49)=\psi'(54)$ and we are forced to define
$\psi''(54)=24+24$ (second block of steps). At this point we have
the equality $\psi''(54)=\psi''(59)$, but we have no more space to
modify $\psi''(59)$.
\end{rem}

The following example shows that the condition of Theorem
\ref{main} is not necessary.

\begin{ex}
Let $S=\langle 16, 17, 35, 71\rangle$. Its Hilbert function $H_R$
is non-decreasing; in fact it assumes the values
$$1,4,8,10,10,11,11,12,12,13,13,14,14,15,15, 16, \rightarrow \ .$$
For $h=3$ we have $|D_h|=|\{52, 70, 88, 106, 142\}|=5 > h+1$,
hence the condition of the theorem is not fulfilled; on the other
hand, computing $C_3$ we get $\{51, 69, 87, 105, 123, 141, 159\}$;
hence $|3M \setminus 4M|-|2M\setminus 3M|=10-8=2=|C_3|-|D_3|$.
\end{ex}

The next result is an immediate consequence of Corollary
\ref{mainb} and Proposition \ref{car}.

\begin{cor}\label{cor0}

If \ $|\left\{\omega_i\in\ap(S)\ :\ a_i>b_i\right\}| \leq 3$,
then $H_R$ is non-decreasing.
\end{cor}

\proof

By Proposition \ref{car}, we have that, if an element $s \in S$
belongs to $D_h$ it must be of the form $s=\omega_i+kg_1 \in
(h-1)M \setminus hM$, for some $\omega_i$ such that $a_i>b_i$ and
$k \in \mathbb{N}$. Furthermore, $\ord(s)=h-1$ implies that, if $k
\neq k'$, then $\omega_i+k'g_1$ cannot have order $h-1$. Thus
$|D_h| \leq 3$ for every $h \in \{2, \dots, r\}$ and from
Corollary \ref{mainb} we get the thesis.
\endproof

We notice that a particular case of the previous corollary is the
well known fact that, if $g_1=4$ (i.e. the multiplicity of $R$
equals $4$), then $H_R$ is non-decreasing.

\medskip

The proof of the main theorem provides a computationally more
efficient condition for $S$, in order to get that $H_R$ is not
decreasing.

\begin{prop}\label{prop}

If $H_R$ is decreasing, then there exists an index $j \geq 2$ such
that $|C_h|\geq h+1$, for every $2 \leq h \leq j$.

\end{prop}

\proof From Corollary \ref{mainb}, there exists an index $h \geq
2$ such that $|D_h|>h+1$. Once we select $h+1$ elements in $D_h$,
we can use the function defined in the proof of Theorem \ref{main}
in order to find $h+1$ different elements in $C_h$. Since each of
these elements in $C_h$ is part of a maximal representation, we
can choose $h$ different maximal representation in $C_{h-1}$ using
the same argument we used in the proof of Theorem \ref{main}.
Again, we can select $h$ elements in $C_{h-1}$ and find $h-1$
elements in $C_{h-2}$, and so on.
\endproof

\begin{rem}
In particular, the index $h$ of the thesis can be chosen as the
index where the Hilbert function decreases. Hence, this could give
a criterion to establish that the Hilbert function is
non-decreasing without computing the cardinalities of $(h-1)M
\setminus hM$ for all the levels $h$. For example, we obtain that,
if $|C_2|< 3$, then $H_R$ is non-decreasing. This fact can be
translated immediately in terms of the Ap\'ery set, as follows.
\end{rem}

\begin{cor}\label{cor}

If $H_R$ is decreasing, then necessarily
\[|\{\omega_i \in\ap(S) :\ b_i=2\}|\geq 3.\]

\end{cor}

\proof By definition $C_2=\{s \in 2M \setminus 3M: s-g_1 \notin
M\setminus 2M,\}$; now, if $s-g_1 \notin M$, $s$ necessarily
belongs to the Ap\'ery set; hence $C_2=\{\omega_i \in\ap(S)\ :\
b_i=2\}$.
\endproof

Notice that Example \ref{es} shows that the condition of
Proposition \ref{prop} is not sufficient, as e.g. $|C_2|=7$, but
$H_R$ is non-decreasing.

The next result is of some interest, since for all the known
examples of decreasing $H_R$ for numerical semigroup rings, the
Hilbert function decreases at first possible step, i.e. for $h=2$
(see, e.g. \cite{MT} and \cite{HW}). Let $e.d.(R)$ be the
embedding dimension of $R$, i.e. the cardinality of the minimal
set of generators for $S$.

\begin{cor}\label{cor2}

If $H_R$ is decreasing at $h=2$, then necessarily
\[e.d.(S)> 5.\]

\end{cor}

\proof

Let $S=\langle g_1, g_2, \dots , g_n\rangle$ (where the generators
are listed as usual in increasing order). If $H_R$ is decreasing
at $h=2$, from Theorem \ref{main}, we get $|D_2|> 3$. Since
$D_2=\{g_j \ : \ g_j+g_1 \in 3M\}$ and since $g_1+g_1$ and
$g_2+g_1$ are in $2M\setminus 3M$, the thesis follows immediately.
\endproof

For a one-dimensional C-M local ring, if $e(R)-e.d.(R) \leq 2$,
the Hilbert function is not decreasing (see e.g. \cite[Theorem
4.8]{Ro}); our last result shows that, in the case of numerical
semigroup rings, for small embedding dimension, this difference
can be increased.

\begin{cor}\label{ed45}
If $e.d.(S)=4,5$ and $g_1 \leq 8$, then $H_R$ is non decreasing.
\end{cor}

\proof By the previous corollary $H_R$ cannot decrease at $h=2$.
If $H_R$ would be decreasing ad $h \geq 3$, then, by Corollary
\ref{cor}, we would get $|C_2| \geq 3$ and $|C_3| \geq 4$. Now
$C_3=\{s \in 3M \setminus 4M : \ s-g_1 \notin 2M\setminus 3M
\}=\{s \in 3M \setminus 4M : \ s-g_1 \in M\setminus 2M\}\cup\{s
\in 3M \setminus 4M : \ s \in \ap(S)\}$.

If $e.d.(S)=4$, we have $\ap(S) \supseteq \{0, g_2, g_3, g_4\}\cup
C_2$; moreover $\{s \in 3M \setminus 4M : \ s-g_1 \in M\setminus
2M\} \subset \{g_3+g_1, g_4+g_1\}$. Hence, to obtain $|C_3| \geq
4$, we need at least two more elements in $\ap(S)$, so $g_1 \geq
9$.

If $e.d.(S)=5$, we have $\ap(S) \supseteq \{0, g_2, g_3, g_4, g_5
\}\cup C_2$; moreover $\{s \in 3M \setminus 4M : \ s-g_1 \in
M\setminus 2M\} \subset \{g_3+g_1, g_4+g_1, g_5+g_1\}$. Hence, to
obtain $|C_3| \geq 4$, we need at least one more elements in
$\ap(S)$, so, again, $g_1 \geq 9$.
\endproof

\begin{footnotesize}
\noindent Corresponding author: Marco D'Anna. e-mail address:
mdanna@dmi.unict.it; fax: 0039 095 330094

\end{footnotesize}


\begin{thebibliography}{99}\footnotesize



\bibitem{A-M-S} F.~Arslan, P.~Mete, M.~Sahin {\em Gluing and Hilbert function of monomial curves},
Proc. Amer. Math. Soc. {\bf 137} (2009), 2225-2232.

\bibitem{BDF} V. Barucci, D. D. E. Dobbs, M. Fontana
{\em Maximality properties in numerical semigroups and applications to one-dimensional
analitically irreducible local domains}, Mem. Amer. Math. Soc., Vol 125, {\bf 598} (1997).

\bibitem{BF} V. Barucci, R. Fr\"oberg,
{\em Associated graded rings of one dimensional analytically irreducible rings}, J.~Algebra {\bf 304} (2006) n.1, 349-358.

\bibitem{CJZ} T. Cortadellas Benitez, R. Jafari, S. Zarzuela Armengou
{\em On the Apery sets of monomial curves} Semigroup Forum {\bf
86} n.2 (2013), pp 289-320

\bibitem{DMM} M.~D'Anna, M.~Mezzasalma, V.~Micale,
{\em On the Buchsbaumness of the associated graded ring of a one-dimensional local ring},
Comm. Alg. {\bf 37} n.5 (2009), 1594 - 1603.

\bibitem{E} J. Elias, {\em The conjecture of Sally on the Hilbert functions for curve singularities}, J. Algebra
{\bf 160} (1993), 42–49.

\bibitem{EB} J. Elias, J. Mart´ýnez-Borruel, {\em Hilbert polynomials and the intersection of ideals}, Contemp.
Math. {\bf 555} (2011), 63–70.

\bibitem{Ga} A. Garcia, {\em Cohen-Macaulayness of the Associated Graded Ring of a Semigroup Ring},
Comm. Algebra {\bf 10} (1982), 393-415.

\bibitem{GR} S. K. Gupta, L. G. Roberts, {\em Cartesian squares and ordinary singularities of curves},
Comm. Algebra {\bf 11} (1983), 127–182.

\bibitem{HW} J. Herzog, R. Waldi, {\em A note on the Hilbert function of a one-dimensional Cohen-Macaulay
ring}, Manuscripta Math. {\bf 16} (1975), no. 3, 251-260.

\bibitem{JZ} R. Jafari, S. Zarzuela Armengou {\em On monomial curves obtained by
gluing}, Semigroup Forum {\bf 88} n.2 (2014), pp 397-416

\bibitem{MT} S.~Molinelli, G.~Tamone {\em On the Hilbert function of certain rings of monomial curves},
Journal of Pure and Applied Algebra
{\bf 101} (1995), 191-206.

\bibitem{N} D.~G.~Northcott, {\em On the notion of a first neighbourhood ring},
Proc.~Camb.~Phil.~Soc. {\bf 53} (1959), 43-56.

\bibitem{NR} D.~G.~Northcott, D.~Rees, {\em Reduction of ideals in local rings},
Proc.~Camb.~Phil.~Soc. {\bf 50} (1954), 145-158.

\bibitem{O} F. Orecchia, {\em One-dimensional local rings with
reduced associated graded ring and their Hilbert functions},
Manuscripta Math. {\bf 32} (1980), 391–405.

\bibitem{P-T} D.P. Patil, G. Tamone, {\em CM defect and Hilbert functions of monomial
curves}, Journal of Pure and Applied Algebra {\bf 215} (2011),
1539-1551.

\bibitem{Ro} M. E. Rossi, {\em Hilbert function of Cohen-Macaulay local rings},
Commutative Algebra and its Connections to Geometry (PASI 2009), Contemporary Mathematics (2010).

\end{thebibliography}
\end{document}